 \documentclass{elsart}
\usepackage{ifpdf}
\usepackage{amsmath,amssymb,latexsym,footnote}
\usepackage{graphicx,amssymb,lineno}
\ifpdf
\usepackage[%
  pdftitle={Instructions for use of the document class
    elsart},%
  pdfauthor={Simon Pepping},%
  pdfsubject={The preprint document class elsart},%
  pdfkeywords={instructions for use, elsart, document class},%
  pdfstartview=FitH,%
  bookmarks=true,%
  bookmarksopen=true,%
  breaklinks=true,%
  colorlinks=true,%
  linkcolor=blue,anchorcolor=blue,%
  citecolor=blue,filecolor=blue,%
  menucolor=blue,pagecolor=blue,%
  urlcolor=blue]{hyperref}
\else
\usepackage[%
  breaklinks=true,%
  colorlinks=true,%
  linkcolor=blue,anchorcolor=blue,%
  citecolor=blue,filecolor=blue,%
  menucolor=blue,pagecolor=blue,%
  urlcolor=blue]{hyperref}
\fi

\makeatletter
\def\elsartstyle{%
    \def\normalsize{\@setfontsize\normalsize\@xiipt{14.5}}
    \def\small{\@setfontsize\small\@xipt{13.6}}
    \let\footnotesize=\small
    \def\large{\@setfontsize\large\@xivpt{18}}
    \def\Large{\@setfontsize\Large\@xviipt{22}}
    \skip\@mpfootins = 18\p@ \@plus 2\p@
    \normalsize
}
\@ifundefined{square}{}{\let\Box\square}
\makeatother

\makeatletter  
\def\@@insvline#1#2{{\setbox0\hbox{\m@th$#1\mathrm I$}  
  \rlap{\m@th$#1 \mkern 5mu  
  \vrule height.95\ht0 depth-.005\ht0 width.09\ht0 $}  
  {\mathrm #2} }}
\def\Q{\mathpalette\@@insvline{Q}}

\makeatother  
  \newtheorem{defi}{Definition}
  \newcommand{\bd}{\begin{defi}} 
  \newcommand{\ed}{\end{defi}}

  \newtheorem{lemm}[defi]{Lemma}  
  \newcommand{\bl}{\begin{lemm}}
  \newcommand{\el}{\end{lemm}} 

  \newtheorem{theo}[defi]{Theorem}
  \newcommand{\bt}{\begin{theo}}
  \newcommand{\et}{\end{theo}}

  \newcommand{\bc}{\begin{cor}}
  \newcommand{\ec}{\end{cor}}

  \newtheorem{pro}[defi]{Proposition}
  \newcommand{\bp}{\begin{pro}}
  \newcommand{\ep}{\end{pro}}
  
  \def\div{\operatorname{div}}
  \makeatletter
  \def\proof{\@ifnextchar[\opargproof{\opargproof[\bf Proof \hfil\\ ]}}
  \def\opargproof[#1]{\par\noindent {\bf #1 }}
  
  \makeatother

\begin{document}
\begin{frontmatter}

%\begin{minipage}{\textwidth}
%\elsartstyle
%
\renewcommand{\thempfootnote}{\fnstar{mpfootnote}}
\leftskip=2pc
\title{Results for a turbulent system with unbounded viscosities:
  weak formulations, existence of solutions, boundedness, smoothness}
%\footnotetext[1]{\upshape Expanded version of a talk presented at the
%Singapore Meeting on Particle Physics (Singapore, August 1990).}
%
\thanks{This work was partially supported by the Swiss National Science Foundation under Grant Number 200020-100051/1.}
\author{P. Dreyfuss}
\address{Institut Elie Cartan, UMR 7502,  
Nancy-Universit\'e, CNRS, INRIA, France}
\ead{Pierre.Dreyfuss@iecn.u-nancy.fr}
\ead[url]{www.iecn.u-nancy.fr/\~{ }dreyfuss }

%\end{center}
%
%\renewcommand{\thempfootnote}{\astsymbol{mpfootnote}}
%\footnotetext[1]{Corresponding author.}
%\setbox0=\hbox{\footnotesize 1}
%\edef\thempfootnote{\hskip\wd0}
%\footnotetext[0]{\textit{Email adresse:} Pierre.Dreyfuss@iecn.u-nancy.fr}
%\footnotetext[1]{\textit{URL:} www.iecn.u-nancy.fr/\~{ }dreyfuss }
%\renewcommand{\thempfootnote}{\arabic{mpfootnote}}
%\footnotetext[1]{

%\bigskip
%\leftskip=0pt

%\hrule\vskip 8pt
%\begin{small}
\begin{abstract}
We consider a circulation system arising in
 turbulence modelling in fluid dynamics with unbounded eddy
 viscosities. Various notions of weak solution are 
considered and compared. We establish existence and regularity
 results. In particular we study the boundedness of weak solutions. We
 also establish an existence result for a classical solution.
\end{abstract} 

%\begin{msc}
%19735
%\end{msc}
%\begin{classification}
%\end{classification}

\end{frontmatter}

%\bigskip
 
\section{Introduction}

Let $\Omega$ be an open bounded set in $\mathbb R^3$, with a Lipschitz boundary. We consider 
the following turbulent circulation model : 
$$
\text{(P)}  \quad\left\{
	\begin{array}{l} 
	-\div(\nu(k)\nabla u) = f \quad \text{in } \Omega \\ 
        -\div(a(k)\nabla k) = \nu(k)| \nabla u|^2 \quad \text{in } \Omega \\ 
 	u = 0 \quad \text{on } \partial \Omega \\ 
	k = 0 \quad \text{on } \partial \Omega \\ 
	\end{array}
	\right. 
$$ 
Here $f, a$ and $\nu$ are given, and the functions $u,k: \ \Omega \to
\mathbb R$ are the unknowns. \\

We study Problem $(P)$ under the following main assumption: 

$$
\text{($H_0$)}  \quad\left\{
	\begin{array}{l}
	  f \in L^r(\Omega), \ \text{with } r > \frac{3}{2} \\ 
	  a, \nu: \mathbb R^{+} \to \mathbb R^{+} \text{ are continuous}\\ 
	  \exists \ \delta > 0: \ a(s),\nu(s) \geq \delta \quad
	  \forall s \in \mathbb R^{+}
	  \end{array}
	\right.
$$
%\end{minipage}

Problem (P) is a simplified scalar version of the RANS model arising
in oceanography (see \cite{lewan,lewan2,num2}): the function $u$ is an 
idealisation of the mean velocity of the fluid and $k$ is the
turbulent kinetic energy. The mathematical analysis of (P) is a step 
towards better understanding the RANS model. Various studies were
made in this direction. Some existence results were established in
\cite{lewan,gal}. \\ 
In this paper we focus on the case where the viscosity functions $a$
and $\nu$ are not a priori bounded. In fact (see \cite{lewan2,gal}), in the
relevant physical situation, we have 
$$
\text{($H_p$)}  \quad\left\{
	\begin{array}{l}
	  a(s) = a_1 + a_2 \sqrt{s} \\ 
	  \nu(s) = \nu_1 + \nu_2 \sqrt{s}
	  \end{array}
	\right.
$$ 

We will establish an existence result for a weak solution for (P) under
less restrictive assumptions than in \cite{gal}. An important feature is that our assumptions are
satisfied under ($H_p$), contrarily to the assumptions made in \cite{gal}. \\ 
Moreover we give additional regularity results for the weak solution
we obtain. In particular, under ($H_0$) and the following additional
assumption: $a$ is proportional to $\nu$, $\partial \Omega$ is of class
$\mathcal{C}^{2,\alpha}$, $f \in
\mathcal{C}^{0,\alpha}(\overline{\Omega})$ and $\nu \in
\mathcal{C}^{1,\alpha}(\mathbb R^{+})$, we prove the existence of a classical solution for
(P). \\ 

We also compare our results with the results presented in
\cite{lewan}. \\ 

Another feature of our work is to considere various notions
of weak solution for Problem (P): $W$-solution, $H$-solution,
distributional solution, renormalized solution, 'energy solution',
classical solution. We give some relations between these notions.

\subsection{Notions of weak solution for (P)}

We can reformulate equation (P).2 by using the Kirchoff
transform. Let 
$$A(s):=\int_0^s a(t)dt.$$
Instead of (P).2, we can consider 
$$(P).2' \qquad -\Delta K = \nu \circ A^{-1}(K)|\nabla u|^2 \quad
\text{on } \Omega,$$
where $K=A(k)$. \\ 
In fact, from every distributional solution $K \in W^1(\Omega)$ of (P).2' we
obtain a distributional solution $k$ of (P).2 by setting
$k=A^{-1}(K)$. This property is related to the facts that $A$ is invertible, $A^{-1}(0)=0$ and 
$|A^{-1}(s)| \leq C.s$ (this can be seen by using the assumptions made
on $\nu$ in ($H_0$)). \\ 

The situation is more complicated for equation (P).1, where the a
priori unbounded coefficient $\nu(k)$ appears in the principal part of
the operator and cannot be removed. Hence we have to restrict $u$ to
satisfy the energy condition 
\begin{equation} 
\int_{\Omega} \nu(k)|\nabla u|^2 < \infty. \label{0.1}
\end{equation}
Nevertheless we will see later on that various non equivalent notions
of weak solution can be considered for (P).1. \\ 
We will introduce the notions of W-solution and H-solution. It is
also possible to consider the notion of renormalized solution
(see \cite{lewan} chap.5). In \cite{gal} the authors defined another notion that they
call energy solution. \\ 
We will give some relations between these notions in the Appendix I. \\ 

Remark now that under the restriction (\ref{0.1}), the right hand side
in (P).2 (or in (P).2') is only a priori in $L^1(\Omega)$. Hence 
(see \cite{BG}) it is natural to seek $k$ in the space 
$\cap_{p<3/2} W_0^{1,p}(\Omega)$. \\ 
We want to find a function $u$ vanishing on $\partial \Omega$ that
satisfies the energy condition (\ref{0.1}). This leads to considering the 
following spaces:
\begin{eqnarray}
W_k &=& \Big\{ v \in H^1_0(\Omega) : [v]_k < \infty \Big\} \notag \\ 
H_k &=& \text{ closure of } \mathcal{C}_c^{\infty}(\Omega) \text{ with
  respect to } [.]_k \notag 
\end{eqnarray}
where we used the notation 
$$[v]_k = \big( \int_{\Omega} \nu(k)|\nabla v|^2 \big)^{1/2}.$$
For any measurable function $k$, the map $[.]_k$ defines a
norm on $W_k$. In the general situation $H_k$ and $W_k$ are not
equal. Moreover $W_k$ is not necessarily complete and a function in 
$H_k$ does not always have a uniquely defined gradient (see \cite{zhikov}). If we 
assume that $\nu(k) \in L^1(\Omega)$ then $W_k$ is complete and in
fact $H_k \subset W_k$ are Hilbert spaces (see \cite{drey,zhikov,cassano}) when they are equipped with
the scalar product 
$$(v,w) = \int_{\Omega} \nu(k) \nabla v \nabla w.$$ 
Consequently, we will consider the following two distinct notions of
solution for (P).1:
$$ 
\begin{array}{l}
u \text{ is called a } H_k\text{-solution of (P).1 if } u \in H_k
\text{ and} \\
\int_{\Omega} \nu(k) \nabla u \nabla v = \int_{\Omega} f v \quad 
\forall v \in H_k
\end{array}$$
$$
\begin{array}{l}
u \text{ is called a } W_k\text{-solution of (P).1 if } u \in W_k
\text{ and} \\
\int_{\Omega} \nu(k) \nabla u \nabla v = \int_{\Omega} f v \quad 
\forall v \in W_k
\end{array}$$ 
Finally, we define the following notions of weak solution for (P):
$$ 
\begin{array}{l}
(u,k) \text{ is called a {\bf H-solution} of (P) if } \\ 
k \in \cap_{p<3/2} W_0^{1,p}(\Omega), \ u \in H_k, \\ 
k \text{ is a distributional solution of (P).2 and } u \text{ is a } H_k
\text{-solution of (P).1} 
\end{array}$$
$$
\begin{array}{l}
(u,k) \text{ is called a {\bf W-solution} of (P) if } \\ 
k \in \cap_{p<3/2} W_0^{1,p}(\Omega), \ u \in W_k \\ 
k \text{ is a distributional solution of (P).2 and } u \text{ is a } W_k
\text{-solution of (P).1} 
\end{array}$$

\subsection{Main results}

Let ($H_1$) and ($H_2$) denote the following conditions:
$$ 
\begin{array}{l}
\text{(}H_1\text{)}\quad \exists \gamma > 0 : \ a(s) \geq \gamma
  \nu(s) \quad \forall s \in \mathbb R^{+} \\ 
\text{(}H_2\text{)}\quad \exists \gamma > 0 : \ a(s) = \gamma
  \nu(s) \quad \forall s \in \mathbb R^{+}
\end{array}
$$

We will establish: 

\begin{theo}\label{theo1}
Assume that ($H_0$) and ($H_1$) hold. Then there exists at least one 
W-solution $(u,k)$ for (P) such that  
\begin{equation}
u\in L^{\infty}(\Omega) \text{ and } 
\int_{\Omega}a(k)|\nabla k|^2 < \infty. \label{reg1}
\end{equation}
\end{theo}

\begin{cor}\label{cor1}
Assume that in addition to ($H_0$) and ($H_1$) we have 
\begin{equation}
\exists \ \nu_0 > 0: \quad \nu(s) \leq \nu_0(1 + s^6), \quad \forall
s \in \mathbb R^{+}. \label{cocor1}
\end{equation}
Then the W-solution $(u,k)$ given in Theorem \ref{theo1} is a 
distributional solution of (P).
\end{cor}

\begin{theo}\label{theo2}
Assume that ($H_0$) and ($H_2$) hold. Then the W-solution $(u,k)$ given in 
Theorem \ref{theo1} satisfies 
\begin{equation} 
u,k \in \mathcal{C}^{0,\alpha}(\overline{\Omega}), \quad 
\text{for some } \alpha \in (0,1). \label{reg2}
\end{equation}
Moreover $(u,k)$ is also a H-solution of (P) (and in fact a classical
weak solution). \\ 
If in addition to ($H_0$) and ($H_2$) we assume that 
$\partial \Omega$ is of class
$\mathcal{C}^{2,\alpha}$, $f \in
\mathcal{C}^{0,\alpha}(\overline{\Omega})$ and $\nu \in
\mathcal{C}^{1,\alpha}(\mathbb R^{+})$ then 
$$u,k \in \mathcal{C}^{2,\beta}(\overline{\Omega}), \quad 
\text{for some } \beta \in (0,1),$$
and $(u,k)$ is a classical solution of (P).
\end{theo}

\subsection{Discussion of the results} \label{sub1}

In Theorem \ref{theo1} we give an existence result of a
W-solution. We next give some regularity results: firstly the 
property (\ref{reg1}) and secondly (in Theorem \ref{theo2}) the property
(\ref{reg2}). Finally, in Theorem \ref{theo2} we give an existence result for a 
classical solution for (P). \\ 

The main previous studies of Problem (P) are presented in 
\cite{lewan} chap. 5 and in \cite{gal}. \\ 

In \cite{lewan} chap.5, the authors prove the existence of a
renormalized solution for (P) under the assumptions
($H_0$) and ($H_2$). It seems that
their proof also works under ($H_0$) and ($H_1$). Nevertheless the
notion of renormalized solution is very weak. A renormalized solution 
$(u,k)$ for (P) is a distributional solution if $\nu(k) \in 
L^{\infty}(\Omega)$, whereas a H- or a W-solution is a distributional
solution if $\nu(k) \in L^1(\Omega)$ (see the Appendix I). \\ 

In \cite{gal} the authors introduced a notion of solution that they
call 'energy solution' (see the Appendix I). In fact an 'energy solution' is a W-solution
which satisfies an additional property ensuring that $H_k = W_k$ 
(the additionnal property imposed is sufficient 
but not necessary to have this equality). Under this point of view an 
'energy solution' is slighty stronger than a W-solution. However, their 
existence result is obtained by assuming complicated conditions on the 
coefficients $a$ and $\nu$ which are not exactly satisfied in the
physically relevant situation ($H_p$), but only in the following approximate
situation: 
$$ 
\text{($H_p'$)}  \quad\left\{
	\begin{array}{l}
	  \text{for some }\epsilon > 0 \text{ we have: } \\ 
	  a(s) = a_1 + a_2 \sqrt{s+\epsilon} \\ 
	  \nu(s) = \nu_1 + \nu_2 \sqrt{s+\epsilon}
	  \end{array}
	\right. 
$$ 
On the contrary, our assumptions in Theorem \ref{theo1} and Corollary \ref{cor1} are
very simple, and they are satisfied in ($H_p$). \\ 
Note also that we establish the regulartity property 
(\ref{reg1}) which are not established in \cite{gal} (or in
\cite{lewan}). \\ 
In the Appendix I we also give a new existence result
for an 'energy solution'. \\ 

In Theorem \ref{theo2} we assume that ($H_0$) and ($H_2$) hold. These
assumptions are fulfilled in the physical situation ($H_p$) if 
$a_2\nu_1 = a_1\nu_2$. We then prove that $u$ and $k$ are 
H\"{o}lder continuous. In particular we give here a positive answer to a
central question put in \cite{gal} : $k$ is bounded. Note that in 
this situation we clearly have $W_k = H_k$. \\ 

We next establish the existence of a classical
solution for Problem (P) by assuming some differentiability properties for 
$a$ and $\nu$. These properties are fulfilled in the situation ($H_p'$) if 
$a_2\nu_1 = a_1\nu_2$. \\ 
It seems that this result is completely new: the existence of a
classical solution for (P) was not studied in any previous work. 

\subsection{Organization of the paper}

In the sequel $n$ will always denote an arbitrary integer
greater or equal to one, and $C$ (possibly with subscript) will denote 
a positive real that does not depend on $n$, but that can differ from
one part to another. \\ 
We always consider the space $H^1_0(\Omega)$ equipped with the
gradient norm. \\ 
The condition ($H_0$) is always assumed.\\§ 

$\bullet \ $In section \ref{sec1} we introduce an approximate sequence $(u_n,k_n)$
of solutions obtained by truncating the coefficients $a$ and $\nu$. \\ 
We immediatly obtain the basic estimates :
\begin{eqnarray*}
\int_{\Omega} \nu_n(k_n)|\nabla u_n|^2 &\leq& C \\ 
\forall p<\frac{3}{2}: \int_{\Omega} \big| a_n(k_n) \nabla k_n \big|^p
&\leq& C 
\end{eqnarray*} 
The point is that we establish the following fundamental estimates:
\begin{eqnarray*} 
\|u_n\|_{L^{\infty}(\Omega)} &\leq& C \\
\int_{\Omega} a_n(k_n)|\nabla k_n|^2 &\leq& C \qquad (*)
\end{eqnarray*}

The first estimate above is proved by developping further a technique due to Stampacchia. \\ 
The second is obtained under the assumption $(H_1)$. The proof is
based on the following idea: if $(u,k)$ is a solution of (P), we
formally have \footnote{We thank Michel Chipot for this remark}
\begin{equation}
\nu(k)| \nabla u|^2 = \underbrace{-\text{ div}(\nu(k) \nabla u).u}_{=fu} + 
\text{ div}(\nu(k) u \nabla u). \label{idee}
\end{equation}
In other words one can hope that the second member in the second
equation in (P) is more regular than it seems. \\
In fact, we prove that a similar relation to (\ref{idee}) holds for
the approximate sequence. By using next that $(u_n)$ is uniformly
bounded in $L^{\infty}(\Omega)$, we obtain $(*)$ which is the key
estimate to prove Theorem \ref{theo1}. \\ 

$\bullet \ $In section \ref{sec2} we extract from $(u_n,k_n)$ a 
subsequence converging to some element denoted by $(u,k)$. Under 
the assumptions ($H_0$) and ($H_1$), we directly obtain that 
$$u \in H^1_0(\Omega) \cap L^{\infty}(\Omega), \quad k \in
H^1_0(\Omega).$$
We prove that moreover we have: 
$$\int_{\Omega} \nu(k)|\nabla u|^2 < \infty, \quad \int_{\Omega}
a(k)|\nabla k|^2 < \infty.$$

$\bullet \ $In section \ref{sec3} we pass to the limit in the approximating
Problems. In a first step we prove that $u$ is a $W_k$-solution of
(P).1. To do this, we use the test functions $v=h_q(k_n)\varphi$ 
(where $\varphi \in W_k \cap L^{\infty}(\Omega)$ and $(h_q)$ is a
sequence of functions that cut off the large values), and we pass to the
limits $n \to \infty$, $q \to \infty$.\\ 
We next prove that the energies of the approximating sequence converge to
the energy $\int_{\Omega} \nu(k)|\nabla u|^2$. \\ 
Finally we can pass to the limit in the second equation in order to
prove that $k$ is a distributional solution of (P).2. We then obtain 
Theorem \ref{theo1} and Corollary \ref{cor1} follows. \\ 

$\bullet \ $In section \ref{sec4} we assume that ($H_0$) and ($H_2$)
hold. In a first step we obtain the estimate 
$$\|k_n\|_{L^{\infty}(\Omega)} \leq C.$$ 
Hence $k \in L^{\infty}$ and by using the De Giorgi-Nash Theorem we
prove the H\"{o}lder continuity of $u$ and $k$. \\ 
Next, by assuming additional regularity on $\nu$, $\partial
\Omega$ and $f$ we can apply the Schauder's estimates and we 
prove Theorem \ref{theo2}. \\ 

$\bullet \ $In the Appendix I we study some relations between the
notions of $W$-solution, $H$-solution,
distributional solution, renormalized solution and 'energy solution'
for Problem (P). We continue the discussion begun in Subsection
\ref{sub1} and we also establish a new existence result for an 'energy
solution' for Problem (P). \\ 
In the Appendix II we recall some basic properties of H\"older continuous functions.

\section{Approximating sequence and estimates} \label{sec1}

We assume that ($H_0$) holds and we set 
\begin{eqnarray}
\nu_n(s) &=& T_n(\nu(s)) \label{1.1} \\ 
a_n(s) &=& T_n(a(s)) \label{1.2},
\end{eqnarray}
where $T_n$ is the truncated function defined by 
$T_n(t)=\text{min}(n,t)$. \\ 
We consider the Problem of finding $(u_n,k_n) \in (H^1_0(\Omega))^2$
such that 
$$
\text{($P_n$)}  \quad\left\{
	\begin{array}{l} 
	\int_{\Omega} \nu_n(k_n)\nabla u_n \nabla v = 
	\int_{\Omega} f v \quad \forall v \in H^1_0(\Omega) \\ 
        \int_{\Omega} a_n(k_n)\nabla k_n \nabla \varphi  =
	\int_{\Omega} T_n\big(\nu_n(k_n)| \nabla u_n|^2 \big) \varphi
	\quad \forall \varphi \in H^1_0(\Omega) 
 	\end{array}
        \right.
$$ 
For any $n\geq 1$, Problem ($P_n$) is well posed because 
$a_n,\nu_n \in L^{\infty}(\mathbb R)$ and $a_n^{-1},\nu_n^{-1} \in
L^{\infty}(\mathbb R)$ by construction. \\ 
It is proved in \cite{gal} that a solution $(u_n,k_n)$ exists for any
$n\geq 1$. Moreover, the following basic properties were established: 
\begin{eqnarray}
k_n &\geq& 0 \label{1.3} \\ 
\int_{\Omega} \nu_n(k_n)|\nabla u_n|^2 &\leq& C_1 \label{1.4} \\ 
\forall p<\frac{3}{2}: \int_{\Omega} \big| a_n(k_n) \nabla k_n \big|^p
&\leq& C_2 \label{1.5} 
\end{eqnarray}
We now establish 
\begin{lemm}\label{lem1}
The sequence $u_n$ is uniformly bounded in the 
$L^{\infty}(\Omega)$-norm, that is, 
\begin{equation} 
\|u_n\|_{L^{\infty}(\Omega)} \leq C_3 \label{1.6}
\end{equation}
\end{lemm}
Before proving this lemma we point out that the assumption 
$f \in L^r(\Omega), \ \text{with } r > \frac{3}{2}$ made in ($H_0$) implies that   
\begin{equation} \label{propf}
f \in W^{-1,\rho}(\Omega), \ \text{with } \rho=\frac{3r}{3-r} > 3.
\end{equation}
This last property is easy to prove by using the Sobolev injection Theorem. 

\begin{proof}
We will obtain the estimate (\ref{1.6}) by using the technique
presented on p.108 in \cite{stam}. \\ 
In order to prove that $C_3$ is independent of $n$ we have to detail
the technique of Stampacchia. \\ 
Let 
$$b_n(u,v):=\int_{\Omega} \nu_n(k_n) \nabla u \nabla v.$$
Recall that $f$ satisfies (\ref{propf}) and then by using a
classical result (see \cite{brezis}) there exists 
$g \in (L^{\rho}(\Omega))^3$ such that $-\div(g)=f$ and 
$\|g\|_{(L^{\rho}(\Omega))^3} \leq C \|f\|_{L^r(\Omega)}$. \\ 
Hence the sequence $u_n$ satisfies 
\begin{equation}
b_n(u_n,v) = \int_{\Omega} g \nabla v \quad \forall v \in
H^1_0(\Omega). \label{lem1.1}
\end{equation}
For $s \geq 0$, we define the measurable set $A_n(s) \subset
\Omega$ by setting 
$$A_n(s) = \big\{ x \in \Omega : \ |u_n(x)| \geq s \big\}.$$
We also introduce 
\begin{equation}
\varphi:=\max\left(|u_n|-s,0\right)sgn(u_n). \label{lem1.2}
\end{equation}
It is proved in \cite{stam} that $\varphi \in H^1_0(\Omega)$ and 
\begin{eqnarray}
\nabla \varphi &=& \nabla u_n \quad \text{in } A_n(s) \notag \\ 
\nabla \varphi &=& 0 \quad \text{in } \Omega \setminus A_n(s) \notag 
\end{eqnarray}
By testing (\ref{lem1.1}) with $v=\varphi$, we obtain 
\begin{equation} 
b_n(\varphi,\varphi)=b_n(u_n,\varphi)=\int_{A_n(s)} g \nabla \varphi.
\label{lem1.3}
\end{equation}
Remark now that assumption $\nu(s)\geq \delta > 0$ in ($H_0$) implies
that $\nu_n(k_n) \geq \min(\delta,1)$. Consequently the bilinear form
$b_n$ is uniformly coercive on $H^1_0(\Omega)$. By using this property
together with the H\"older inequality, we obtain from (\ref{lem1.3}): 
$$\|\varphi\|_{H^1_0(\Omega)}^2 \leq \tilde{C}\big( 
\int_{A_n(s)} |g|^2\big)^{1/2} \|\varphi\|_{H^1_0(\Omega)}.$$
Hence by using the Cauchy inequality together with the H\"older inequality we obtain 
\begin{equation}
\|\varphi\|^2_{H^1_0(\Omega)^2} \leq \tilde{C_1}  
\|g\|_{L^\rho(\Omega)}^2 \big| A_n(s) \big|^{\frac{\rho-2}{\rho}}. \label{lem1.4} 
\end{equation}
On the other hand, the Poincar\'e-Sobolev inequality gives 
\begin{equation}
\big( \int_{A_n(s)} |\varphi|^6\big)^{1/3} \leq \tilde{C_2} 
\|\varphi\|_{H^1_0(\Omega)}^2. \label{lem1.5}
\end{equation}
Let now $t > s$. It is clear that $A_n(t) \subset A_n(s)$ and
consequently 
\begin{equation}
\big( \int_{A_n(s)} |\varphi|^6\big)^{1/3} \geq 
\big( \int_{A_n(t)} |\varphi|^6\big)^{1/3} \geq 
\big( \int_{A_n(t)} |t-s|^6\big)^{1/3} \geq 
|t-s|^2 \big|A_n(t)\big|^{1/3} \label{lem1.6}
\end{equation}
We set 
$$\psi_n(s) := \big|A_n(s)\big|, \quad \forall s \geq 0$$
For fixed $n$, $\psi_n$ is a decreasing function, and from the
estimates (\ref{lem1.4})-(\ref{lem1.6}), we obtain 
$$\psi_n(t) \leq \tilde{C_3}|\psi_n(s)|^{\beta}(t-s)^{-6} \quad 
\forall t > s \geq 0,$$
where we have used the notation $\beta:=\frac{3(\rho-2)}{\rho} > 1$ and where  $\tilde{C_3}=\tilde{C_3}(\tilde{C_1},\tilde{C_2},\|f\|_{L^r})$. Both quantity $\beta$ and $\tilde{C_3}$ do not depend on $n$. 
Hence by using Lemma 4.1 in \cite{stam} it follows: 
$$\psi_n(\theta) = 0,$$
where $\theta=2^{\beta/(\beta -1)}\big(\tilde{C_3}|\Omega|^{\beta-1}\big)^{1/6} < \infty$ 
does not depend on $n$. \\ 
This property tells precisely that (\ref{1.6}) holds true with 
$C_3=\theta$.
\end{proof}

Notice that the bilinear form 
$$(u,v) \to \int_{\Omega} a_n(k_n) \nabla u \nabla v,$$
is also uniformly coercive on $H^1_0(\Omega)$. Moreover, the sequence 
$$h_n := T_n\big(\nu_n(k_n)| \nabla u_n|^2 \big)$$ 
is imbedded in $L^{\infty}(\Omega)$. We can then apply again the 
technique of Stampacchia detailed in the proof of lemma \ref{lem1},
and obtain: 
\begin{equation}
\text{for } n \geq 1: \quad k_n \in L^{\infty}(\Omega) \label{1.7}
\end{equation}
Nevertheless the control we have on $\{h_n\}$ is obtained from 
(\ref{1.4}), which gives a uniform bound in the $L^1$-norm for 
the sequence. This is not enough to obtain a uniform estimate for 
$\{k_n\}$ in the $L^{\infty}$-norm. \\ 
However we can establish: 

\begin{lemm}\label{lem2}  
Assume that ($H_0$) and ($H_1$) hold. Then we have 
\begin{eqnarray}
a_n(s) &\geq& \gamma_1 \nu_n(s), \quad \gamma_1=\min(1,\gamma)
\label{1.8'} \\ 
\int_{\Omega} a_n(k_n)|\nabla k_n|^2 &\leq& C_5 \label{1.8}
\end{eqnarray}
\end{lemm}

\begin{proof}
The estimate (\ref{1.8'}) is easy to obtain. Its verification is left
to the reader. \\ 
Let $(u_n,k_n)$ be the chosen approximating sequence. We have from 
(\ref{1.6}) and (\ref{1.7}) that 
$$\forall n \geq 1: \quad u_n,k_n \in H^1_0(\Omega) \cap
L^{\infty}(\Omega)$$
It follows (see \cite{brezis}) that 
$v:=u_n.k_n \in  H^1_0(\Omega) \cap L^{\infty}(\Omega)$ is admissible
for ($P_n$).1 and we get\footnote{more generally: 
$\nu_n(k_n)|\nabla u_n|^2 = \underbrace{-\text{div}(\nu_n(k_n) \nabla
    u_n).u_n}_{=f u_n} + 
\text{div}(\nu_n(k_n) u_n \nabla u_n)$ in $\mathcal{D}'(\Omega)$}
\begin{equation}
\int_{\Omega} \nu_n(k_n)|\nabla u_n |^2 k_n = 
	\int_{\Omega} f u_n k_n -  \int_{\Omega} \nu_n(k_n)u_n\nabla
	u_n \nabla k_n \label{lem2.1}
\end{equation}
By testing ($P_n$).2 with $\varphi=k_n$, we obtain: 
\begin{equation}
\int_{\Omega} a_n(k_n)|\nabla k_n |^2 = 
\int_{\Omega} T_n\big(\nu_n(k_n)| \nabla u_n|^2 \big) k_n \leq 
\int_{\Omega} \nu_n(k_n)|\nabla u_n |^2k_n, \label{lem2.2}
\end{equation}
by using the properties $T_n(s) \leq s$ and (\ref{1.3}). \\ 
Hence, by combining (\ref{lem2.1}) with (\ref{lem2.2}) we have: 
\begin{equation}
I:=\int_{\Omega} a_n(k_n)|\nabla k_n |^2 \leq 
\underbrace{\int_{\Omega} |f u_n k_n|}_{:=II} + 
\underbrace{\int_{\Omega} \big| \nu_n(k_n)u_n\nabla u_n \nabla
  k_n\big|}_{:=III} \label{lem2.3}
\end{equation}
We can estimate the term II as follows: 
\begin{eqnarray} 
II &\leq& C_3 \int_{\Omega} |f k_n| \overset{\text{H\" older Ineq.}}{\leq} C_3 \|f\|_{L^{3/2}}\|k_n\|_{L^3}\notag \\ 
&\overset{\text{Poincar\'e-Sobolev Ineq.}}{\leq}& \tilde{C_1} \|f\|_{L^{3/2}} \big( 
\int_{\Omega} |\nabla k_n |^2\big)^{1/2} \leq  
\frac{\tilde{C_1}}{\delta}\|f\|_{L^{3/2}} \big(\int_{\Omega}
a_n(k_n)|\nabla k_n |^2\big)^{1/2} \notag \\ 
&\overset{\text{Young Ineq.}}{\leq}& \frac{\tilde{C_1}}{\delta} \big(\frac{1}{\epsilon} \|f\|_{L^{3/2}}^2 + \epsilon \int_{\Omega} a_n(k_n)|\nabla k_n |^2 \big) 
\quad \text{for any }\epsilon > 0 \text{ given} \notag \\  
&\leq& \frac{1}{3} \int_{\Omega} a_n(k_n)|\nabla k_n |^2 + 
\tilde{C_2} \|f\|_{L^{3/2}}^2 \notag
\end{eqnarray}
where $\delta > 0$ is the constant given in ($H_0$). The
last inequality was obtained by choosing
$\epsilon=\delta/(3\tilde{C_1})$, using the estimate (\ref{1.6})
and by setting $\tilde{C_2}=3\tilde{C_1}^2/\delta^2$. \\ 

We next estimate the term III:
\begin{eqnarray}
III &=& \int_{\Omega} \big| u_n \sqrt{\nu_n(k_n)} \nabla u_n 
\sqrt{\nu_n(k_n)} \nabla k_n \big| \notag \\ 
&\leq& \tilde{C_3} \int_{\Omega} \big| \sqrt{\nu_n(k_n)} \nabla u_n 
\sqrt{a_n(k_n)} \nabla k_n \big|, \quad 
\tilde{C_3}=C_3 \gamma_1^{-1/2} \notag   \\ 
&\leq& \frac{1}{3}\int_{\Omega} a_n(k_n)|\nabla k_n |^2 + 
\tilde{C_4} \int_{\Omega} \nu_n(k_n)|\nabla u_n |^2, \quad 
\tilde{C_4}=\tilde{C_4}(\tilde{C_3}) \notag
\end{eqnarray}
where $C_3, \gamma_1$ are the constants that appear in (\ref{1.6}) and (\ref{1.8'}). The last
inequality follows from the Young inequality. \\ 
Recall now the inequality (\ref{lem2.3}) and use the estimates
established for the terms II and III. We obtain: 
\begin{equation}
\frac{1}{3}\int_{\Omega} a_n(k_n)|\nabla k_n |^2 \leq 
\tilde{C_2} \|f\|_{L^{3/2}(\Omega)}^2 + \tilde{C_4} \int_{\Omega} \nu_n(k_n)|\nabla u_n |^2. 
\label{lem2.4}
\end{equation}
By using (\ref{lem2.4}) together with (\ref{1.4}) we finally obtain 
(\ref{1.8}).
\end{proof}

\section{Basic convergence results for $(u_n,k_n)$} \label{sec2}

The estimates established in the previous section allow us to extract a
converging subsequence from $(u_n,k_n)$. We have 

\begin{lemm} \label{lem3} \noindent \\  
1. Assume that ($H_0$) holds. Then we can extract a subsequence 
(still denoted by $(u_n,k_n)$) such that 
\begin{eqnarray}
a_n(k_n) \nabla k_n &\rightharpoonup& a(k)\nabla k \ \text{ in }
L^p(\Omega), \ p < \frac{3}{2} \label{2.1} \\
k_n &\rightarrow& k \ \text{ a.e in }\Omega \label{2.2} \\
u_n &\rightharpoonup& u \ \text{ in } H_0^1(\Omega) \label{2.3} \\
u_n &\overset{*}{\rightharpoonup}& u \ \text{ in }
L^{\infty}(\Omega)\label{2.4}
\end{eqnarray}

2. If in addition the condition ($H_1$) is fulfilled then we may assume
   that 
\begin{equation}
k_n \rightharpoonup k \ \text{ in } H_0^1(\Omega) \label{2.5}
\end{equation}
\end{lemm}

\begin{proof}
1. The properties (\ref{2.1}) and (\ref{2.2}) are obtained from
   (\ref{1.5}). The property (\ref{2.3}) is obtained by using the
   estimate (\ref{1.4}) together with the assumption $\nu(s) \geq
   \delta > 0$ in ($H_0$). We establish (\ref{2.4}) from the estimate 
(\ref{1.6}). \\ 

2. By using Lemma \ref{lem2} together with the assumption $a(s)\geq
   \delta > 0$ in ($H_0$) we obtain (\ref{2.5}). Notice that the $k$
   appearing in (\ref{2.1}), (\ref{2.2}) and (\ref{2.5}) is necessarily
   the same in the three situations. 
\end{proof}

We are able to prove additional regularity results for the 
element $(u,k)$ introduced in Lemma \ref{lem3}. For technical reasons
we introduce the sequence $\{h_q\}_{q \in \mathbb N}$ of real functions
defined in \cite{lewan} p. 185. It satisfies: 
\begin{eqnarray}
|h_q(s)| &\leq& 1 \quad \forall (q,s) \in \mathbb N \times \mathbb R
 \label{p1} \\
h_q(s) &=& 0 \quad \text{when } |s| > 2q \label{p0} \\ 
|h_q'(s)| &\leq& \frac{1}{q} \quad \forall q \in \mathbb N, \text{ and
 a.e } s \in \mathbb R
 \label{p2} \\ 
h_q &\underset{q \to \infty}{\rightarrow}& 1\quad \text{uniformly on
  the compacts} \label{p3}
\end{eqnarray}

\begin{lemm}\label{lem4} \noindent \\  
1. Assume that ($H_0$) holds. Then the element $(u,k)$ given in
   Lemma \ref{lem3} satisfies
\begin{equation}
\int_{\Omega} \nu(k)|\nabla u |^2 < \infty \label{2.6}
\end{equation}
2. Assume that in addition ($H_1$) holds. Then 
\begin{equation}
\int_{\Omega} a(k)|\nabla k |^2 < \infty \label{2.7}
\end{equation}
\end{lemm}

\begin{proof}
1. We take over the arguments presented in \cite{lewan} p. 192. \\ 
For $q\geq 1$, we set 
$$\eta_{n,q}:=\big(h_q(k_n)\nu_n(k_n)\big)^{1/2}\nabla u_n$$
Let now $q$ be fixed. The sequence 
$\big\{\big(h_q(k_n)\nu_n(k_n)\big)^{1/2}\big\}_{n\geq 1}$ is uniformly
bounded in $L^{\infty}(\Omega)$. Consequently, $\{\eta_{n,q}\}_{n\geq 1}$
is bounded in $(L^2(\Omega))^3$ and we can extract a subsequence 
weakly convergent to some $\eta_q \in (L^2(\Omega))^3$. \\ 
On the other hand, we have 
\begin{eqnarray}
\big(h_q(k_n)\nu_n(k_n)\big)^{1/2} &\rightarrow&
\big(h_q(k)\nu_n(k)\big)^{1/2} \quad \text{a.e in } \Omega \notag \\
\nabla u_n &\rightharpoonup& \nabla u \quad \text{in } L^2(\Omega),
\notag
\end{eqnarray}
and thus $\eta_q=\big(h_q(k)\nu(k)\big)^{1/2} \nabla u$. \\ 
We now use a classical property of the weak convergence in
$L^2(\Omega)$:
\begin{eqnarray}
\|\eta_q\|_{L^2(\Omega)} &\leq& \liminf\limits_{n\to \infty}
\|\eta_{n,q}\|_{L^2(\Omega)} \leq \liminf\limits_{n\to \infty} 
\big(\int_{\Omega} \nu_n(k_n)|\nabla u_n|^2\big)^{1/2} \notag 
\leq C_1^{1/2}, \notag
\end{eqnarray}
where $C_1$ is a constant independent of $q$ given in (\ref{1.4}). \\ 
By using properties (\ref{p3}) and (\ref{p1}) we can see that 
\begin{eqnarray}
&\eta_q^2& \underset{q\to\infty}{\rightarrow}\nu(k)|\nabla u|^2 \quad
\text{a.e. in }\Omega \notag \\
&\eta_q^2&\leq \nu(k)|\nabla u|^2 \notag
\end{eqnarray}
Hence by the Fatou Lemma we finally obtain: 
$$\int_{\Omega} \nu(k)|\nabla u |^2 \leq \liminf\limits_{q \to \infty}
\|\eta_q\|_{L^2}^2 \leq C_1$$

2. If the additional assumption ($H_1$) holds, then we have the
   estimate (\ref{1.8}) and the previous reasoning allows us to obtain
   (\ref{2.7})
\end{proof}

\section{The proof of theorem \ref{theo1}} \label{sec3}

In the previous section we have proved that under ($H_0$) we can
extract a converging subsequence of $(u_n,k_n)$. If moreover ($H_1$)
holds then the limit $(u,k)$ obtained satisfies: 
\begin{eqnarray}
u &\in& W_k \cap L^{\infty}(\Omega) \label{3.1.a} \\ 
k &\in& H^1_0(\Omega) \quad \text{(and in fact } k\in W_k\text{)} \label{3.1.b}
\end{eqnarray}

\subsection{Passing to the limit in ($P_n$).1}

We recall that the space $W_k$ was defined by 
$$W_k = \Big\{ v \in H^1_0(\Omega) : [v]_k < \infty \Big\}$$
We now establish: 

\begin{lemm}\label{lem5}
Assume that ($H_0$) and ($H_1$) hold. Then the element $(u,k)$
given in Lemma \ref{lem3} satisfies (\ref{3.1.a}), (\ref{3.1.b}) and:
\begin{equation}
\int_{\Omega} \nu(k) \nabla u \nabla v = \int_{\Omega} f v \quad 
\forall v \in W_k \label{3.2}
\end{equation}
\end{lemm}
\begin{proof}
Let $n\geq1, q \in \mathbb N$ and $\varphi \in W_k\cap
L^{\infty}(\Omega)$. We consider the function $v:=h_q(k_n)\varphi$. By
recalling the properties (\ref{p1})-(\ref{p3}) of $h_q$, we can verify
that $h_q(k_n) \in H^1_0(\Omega)\cap L^{\infty}(\Omega)$. Consequently 
$v\in  H^1_0(\Omega)\cap L^{\infty}(\Omega)$. By testing ($P_n$).1
with $v$, we obtain:
\begin{equation}
I:=\int_{\Omega}\nu_n(k_n)h_q(k_n)\nabla u_n \nabla \varphi + 
\underbrace{\int_{\Omega} h_q'(k_n)\nu_n(k_n)\nabla u_n \nabla k_n
  \varphi
}_{:=II} = \underbrace{\int_{\Omega}f h_q{k_n}\varphi}_{:=III}
\label{lem5.1}
\end{equation}
In a first step we fix $q$ and we study the behaviour of terms I,
II and III when $n$ tends to infinity. \\ 
By using the property (\ref{p0}) we see that 
$$|\nu_n(k_n)h_q(k_n)| \leq \max\limits_{s\in[0,2q]} \nu(s) := C_q,$$ 
and by using (\ref{p0}) together with (\ref{2.2}) we obtain 
$$\nu_n(k_n)h_q(k_n) \underset{n\to\infty}{\rightarrow}
\nu(k)h_q(k) \quad \text{a.e in }\Omega.$$
Consequently 
$$\nu_n(k_n)h_q(k_n)\nabla \varphi \underset{n\to\infty}{\rightarrow}
\nu(k)h_q(k)\nabla \varphi \quad \text{in }(L^2(\Omega))^2,$$
and by also employing (\ref{2.3}) we get: 
\begin{equation}
I \underset{n\to\infty}{\rightarrow} \int_{\Omega}\nu(k)h_q(k)\nabla u
\nabla \varphi \label{lem5.2}
\end{equation}
We now estimate II. From (\ref{p2}) we obtain: 
\begin{eqnarray}
II &\leq& \frac{1}{q} \int_{\{q\leq k_n\leq 2q\}} 
\big|\nu_n(k_n) \nabla u_n \nabla k_n \varphi \big| \notag \\ 
&\leq& \|\varphi\|_{L^{\infty}}\frac{C}{q} 
\big( \int_{\Omega} \nu_n(k_n)|\nabla u_n |^2 \big)^{1/2}
\big( \int_{\Omega} a_n(k_n)|\nabla k_n |^2 \big)^{1/2} 
\leq \frac{C}{q}, \label{lem5.3}
\end{eqnarray}
where the second inequality is obtained by using (\ref{1.8'}). \\ 
For the last term we get 
\begin{equation}
III \underset{n\to\infty}{\rightarrow} \int_{\Omega}f h_q(k)\varphi \label{lem5.4}
\end{equation}
By using the estimates (\ref{lem5.2})-(\ref{lem5.4}) together with (\ref{lem5.1})
we obtain that for any fixed $\varphi \in W_k \cap L^{\infty}(\Omega)$
the following holds 
\begin{equation}
\underbrace{\int_{\Omega}\nu(k)h_q(k)\nabla u\nabla \varphi
}_{:=J_1} = 
\underbrace{\int_{\Omega} f h_q(k)\varphi
}_{:=J_2} + \mathcal{O}(\frac{1}{q}). \label{lem5.5}
\end{equation}
We next remark that the integrand in $J_1$ converges for a.e. $x \in
\Omega$ to $\nu(k)\nabla u \nabla \varphi$ when $q$ tends to
infinity. Moreover by using (\ref{p1}) together with the fact that
$\varphi\in W_k$ we can see that the integrand in $J_1$ is dominated
by $|\nu(k)\nabla u \nabla \varphi| \in L^1(\Omega)$. Consequently, by
the Dominated Convergence Theorem we get 
$$J_1 \underset{q\to\infty}{\rightarrow}\int_{\Omega} \nu(k)\nabla u
\nabla \varphi$$
Similarly we can see that 
$$J_2 \underset{q\to\infty}{\rightarrow} \int_{\Omega} f \varphi$$
At this stage we have proved that 
\begin{equation}
\int_{\Omega} \nu(k)\nabla u§\nabla \varphi =  \int_{\Omega} f \varphi
\quad \forall \varphi \in W_k\cap L^{\infty}(\Omega), \label{lem5.6}
\end{equation}
and it remains to show that the condition $\varphi\in
L^{\infty}(\Omega)$ is not necessary. \\ 
Let $\varphi \in W_k$ and $i\in \mathbb N$. We consider 
$\varphi_i \in  W_k\cap L^{\infty}(\Omega)$ given by 
$\varphi_i = T_i(\varphi)$. By using some basic properties of $T_i$
(see \cite{gal}), we see that $|\varphi_i| \leq |\varphi|$, $|\nabla \varphi_i|
\leq |\nabla \varphi|$, $\varphi_i \to \varphi$ a.e
, and $\nabla \varphi_i \to \nabla \varphi$ a.e in
$\Omega$. Consequently, if we take $\varphi_i$ as test function in
(\ref{lem5.6}), we can pass to the limit $i\to\infty$ and we obtain   
(\ref{3.2}).
\end{proof}

In Lemma \ref{lem5} we have showed that $u$ is a $W_k$-solution of
(P).1. In order to prove Theorem \ref{theo1} we have to prove that $k$
is a distributional solution of (P).2. We need first to establish: 

\begin{lemm}\label{lem6}
Assume that ($H_0$) and ($H_1$) hold. Then, in addition to the results
presented in Lemma \ref{lem3}, we may assume: 
\begin{equation}
\nu_n(k_n)|\nabla u_n|^2 \underset{n\to\infty}{\rightarrow} 
\nu(k)|\nabla u|^2 \quad \text{in } L^1(\Omega) \label{3.3}
\end{equation}
\end{lemm}
\begin{proof}
We test ($P_n$).1 with the function $u_n$. By using (\ref{2.3}) we
obtain: 
\begin{equation}
\int_{\Omega}\nu_n(k_n)|\nabla u_n|^2
 \underset{n\to\infty}{\rightarrow} 
\int_{\Omega} f u = \int_{\Omega} \nu(k)|\nabla u|^2, \label{lem6.1}
\end{equation}
where the latter equality is obtained by testing (\ref{3.2}) with
$u$. \\ 
We set $\eta_n:=\sqrt{\nu_n(k_n)}\nabla u_n
$ and $\eta:=\sqrt{\nu(k)}\nabla u$. The relation (\ref{lem6.1}) tells
us that 
\begin{equation}
\|\eta_n\|_{L^2(\Omega)} \underset{n \to \infty}{\rightarrow} \|\eta\|_{L^2(\Omega)}
\label{lem6.2}
\end{equation}
We can next take over the arguments presented in \cite{lewan} Lemma
5.3.4 in order to obtain:
\begin{equation}
\eta_n \underset{n \to \infty}{\rightharpoonup} \eta \quad \text{in } (L^2(\Omega))^2
\label{lem6.3}
\end{equation}
Finally properties (\ref{lem6.3}) and (\ref{lem6.2}) imply that 
the convergence is strong in (\ref{lem6.3}), and (\ref{3.3}) follows. 
\end{proof}

\subsection{The Proofs of Theorem \ref{theo1} and Corollary \ref{cor1}}
 
Assume that ($H_0$) and ($H_1$) hold. In Lemma \ref{lem3} we have
extracted a subsequence $(u_n,k_n)$ which converges in a certain 
sense to an element $(u,k)$. This element has the properties
(\ref{3.1.a})-(\ref{3.1.b}). 
Next we have established (\ref{3.2}). \\ 
Let now $\varphi \in \mathcal{C}_{c}^{\infty}(\Omega)$. By using
(\ref{2.1}) we get: 
\begin{equation} 
\int_{\Omega}a_n(k_n)\nabla k_n \nabla \varphi
\underset{n\to\infty}{\rightarrow} 
\int_{\Omega}a(k)\nabla k \nabla \varphi
 \label{toto1}
\end{equation}
We next remark that the property (\ref{3.3}) ensures that 
\begin{equation}
\int_{\Omega}T_n\big(\nu_n(k_n)|\nabla u_n|^2\big)\varphi 
 \underset{n\to\infty}{\rightarrow} 
\int_{\Omega}\nu(k)|\nabla u|^2\varphi \label{toto2}
\end{equation}

Recall that the sequence $(u_n,k_n)$ satisfies ($P_n$).2. Then
relation (\ref{toto1}) together with (\ref{toto2}) allows to take the
limit in ($P_n$).2. We get: 
\begin{equation}
\int_{\Omega} a(k)\nabla k \nabla \varphi = 
\int_{\Omega} \nu(k)|\nabla u|^2 \varphi \quad \forall 
\varphi \in \mathcal{C}_{c}^{\infty}(\Omega) \label{3.4}
\end{equation}
Thus (P).2 is fulfilled in the distributional sense. \\ 
At this point we have obtained (\ref{3.1.a}), (\ref{3.1.b}), 
(\ref{3.2}) and (\ref{3.4}). The proof of Theorem \ref{theo1} is
complete. \\

Assume now that the condition (\ref{cocor1}) in Corollary \ref{cor1}
is fulfilled. By using (\ref{3.1.b}) together with the Sobolev
Injection Theorem we get $k\in L^{6}(\Omega)$ and thus 
$\nu(k)\in L^{1}(\Omega)$. Then we can conclude the proof of 
Corollary \ref{cor1} by using Proposition
\ref{pro1} in the Appendix I: $(u,k)$ is a distributional solution of (P).    

\section{The proof of Theorem \ref{theo2}} \label{sec4}

We assume in this section that ($H_0$) and ($H_2$) hold. \\ 
In this situation all the results presented in section 2 and section 3
are valid. For technical reasons we slightly modify the definition of
$a_n$ by setting 
\begin{equation}
a_n(s):=\gamma \nu_n(s), \label{titi1} 
\end{equation}
where $\gamma > 0$ is the constant appearing in ($H_2$) and 
$\nu_n$ is defined as before. \\ 
We will now consider Problems ($P_n$) modified by the new definition
(\ref{titi1}) of $a_n$. Nevertheless the modification is very sligth,
and all the results presented in the previous section can be recovered
easly. The verifications are left to the reader. \\ 

We now prove that we have the new estimate:
\begin{equation}
\|k_n\|_{L^{\infty}(\Omega)} \leq C_6 \label{4.1}
\end{equation}
In order to prove this result we set 
\begin{equation} 
\chi_n := k_n + \frac{\gamma}{2} u_n^2, \label{titi2}
\end{equation}
and we remark that ($P_n$).2 leads to 
$$\int_{\Omega} a_n(k_n) \nabla \chi_n \nabla \varphi = 
\int_{\Omega} f u_n \varphi \quad \forall \varphi \in H^1_0(\Omega).$$
Recall that $a_n(k_n)\geq \gamma \text{min}(1,\delta)>0$, $a_n(k_n)\in
L^{\infty}(\Omega)$ and note that the sequence $fu_n$ is uniformly
bounded in $L^r(\Omega)$ with  $r>3/2$. These properties are sufficient (see
the proof of Lemma \ref{lem1}) to get the estimate 
\begin{equation}
\|\chi_n\|_{L^{\infty}(\Omega)} \leq C, \label{titi3}
\end{equation}
where $C$ does not depend on $n$. \\ 
The estimate (\ref{4.1}) is finally obtained by using Lemma \ref{lem1} together
with (\ref{titi3}). \\ 
Consequently, in addition to the properties in Lemma \ref{lem3} we may
assume that 
\begin{equation}
k_n \overset{*}{\rightharpoonup} k \quad \text{in } L^{\infty}(\Omega). 
\label{4.2}
\end{equation}

We will now prove that 
\begin{equation}
u,k \in \mathcal{C}^{0,\alpha}(\overline{\Omega}) \quad \text{for some
} \alpha \in (0,1). \label{4.4}
\end{equation}  

Let $\lambda:=\nu(k)$. We have $\lambda, \lambda^{-1} \in
L^{\infty}(\Omega)$ and 
\begin{equation}
\int_{\Omega} \lambda \nabla u \nabla \phi = \int_{\Omega} f \phi
\quad \forall \phi \in H_0^1(\Omega). \label{tata1}
\end{equation}
Recall also that $f$ have the property (\ref{propf}). Hence 
we can apply the De Giorgi-Nash Theorem (see for instance
\cite{DL} Prop. 6 p.683 or \cite{gil} Th. 8.22 and Th. 8.29). We obtain that 
$u \in \mathcal{C}^{0,\alpha_1}(\overline{\Omega})$ for some $\alpha_1
\in (0,1)$. We next set $\chi:=k+(\gamma/2)u^2$. Then 
$\chi \in H_0^1(\Omega)$ and we have  
\begin{equation}
\int_{\Omega} \frac{\lambda}{\gamma} \nabla \chi \nabla \phi = 
\int_{\Omega} fu \phi \quad \forall \phi \in H_0^1(\Omega).
\label{tata2}
\end{equation}
By using the fact that $u \in L^{\infty}(\Omega)$ in 
(\ref{tata2}), we can again apply the De Giorgi-Nash Theorem to get 
$\chi \in \mathcal{C}^{0,\alpha_2}(\overline{\Omega})$ for some $\alpha_2
\in (0,1)$. Hence also $k$ is H\"older continuous, and (\ref{4.4})
follows. \\ 

Let $\alpha\in (0,1)$ be a generic parameter that can differ from one
part to another. We 
assume now that $\partial \Omega$ is of class
$\mathcal{C}^{2,\alpha}$, $f \in
\mathcal{C}^{0,\alpha}(\overline{\Omega})$ and $\nu \in
\mathcal{C}^{1,\alpha}(\mathbb R^{+})$. \\ 
We will prove the second part of Theorem \ref{theo2} by iterating the Schauder estimates. \\ 
We have $\lambda=\nu(k) \in \mathcal{C}^{0,\alpha}(\overline{\Omega})$
 (see the Appendix II) and then, by applying the Schauder estimate (see \cite{chen} Theorem 2.7
p. 154) on (\ref{tata1}) we get 
$u\in \mathcal{C}^{1,\alpha}(\overline{\Omega})$. Similarily, from
equation (\ref{tata2}) we obtain 
$\chi \in \mathcal{C}^{1,\alpha}(\overline{\Omega})$ 
and thus 
$k \in \mathcal{C}^{1,\alpha}(\overline{\Omega})$. \\ 
Hence (see Appendix II) $\lambda \in  \mathcal{C}^{1,\alpha}(\overline{\Omega})$. By 
iterating again the Schauder estimates (see now Theorem 2.8 p.154 in \cite{chen}) we obtain that 
$u$ and $k$ are in $\mathcal{C}^{2,\alpha}(\overline{\Omega})$. \\ 
Finally we see that $(u,k)$ is a classical solution of (P). Theorem
\ref{theo2} is proven.

\section*{Appendix I: Some relations between the notions of weak solution}

We give here some relations between the various notions of weak
solution: $W$-solution, $H$-solution, distributional solution,
renormalized solution, 'energy solution'.

\subsection*{Comparison with renormalized solution}

We have: 

\begin{pro}\label{pro1} \noindent \\ 
1. Any $W$- or $H$-solution $(u,k)$ of Problem (P) that satisfies in addition  
   $k\in H^1_0(\Omega)$, is also a renormalized solution. \\ 
2. If $\nu(k) \in L^1(\Omega)$ then any $W$- or $H$-solution of
   Problem (P) is a also a distributional solution of (P).
\end{pro}
\begin{proof}
1. Let $(u,k)$ be a $W$-solution of (P). Then the conditions
   (5.2.1)-(5.2.5) in \cite{lewan} chap.5 are satisfied. We have to
prove that (5.2.6) holds. \\ 
Let $h \in \mathcal{C}_{c}^{\infty}(\mathbb R)$ and 
$\phi \in \mathcal{C}_{c}^{\infty}(\Omega)$ be arbitrarily chosen. We
set $v:=h(k)\phi$. Then  $v\in L^{\infty}(\Omega)$ and 
$\nabla v = h(k)\nabla \phi + h'(k)\nabla k \phi$. Let $M < \infty$ be 
such that the support of $h$ being included in $[-M,M]$. We have 
\begin{eqnarray}
\int_{\Omega} \nu(k) h^2(k) |\nabla \phi|^2 &\leq& 
\max\limits_{[0,M]}\nu \|h\|_{L^{\infty}}^2 \int_{\Omega}|\nabla \phi|^2
< \infty \notag \\ 
\int_{\Omega} \nu(k) (h'(k))^2 |\nabla k|^2 |\phi|^2 &\leq& 
\max\limits_{[0,M]}\nu \|h'\|_{L^{\infty}}^2 \|\phi\|_{L^{\infty}}^2\int_{\Omega}|\nabla k|^2
< \infty \notag  
\end{eqnarray}
Hence $v\in W_k$. By testing (\ref{3.2}) with $v$ we obtain the
relation (5.2.6).a in \cite{lewan}. \\ 
We remark that $v$ is also admissible in (\ref{3.4}). This allows us to
obtain the condition (5.2.6).b in \cite{lewan}. Consequently 
$(u,k)$ is a renormalized solution of (P). \\ 
If we consider a $H$-solution $(u,k)$ of (P) we can take over the
previous argument because the function $v$ is now in $H_k$. \\ 
2. If $\nu(k)\in L^1(\Omega)$ then we have 
$\mathcal{C}_{c}^{\infty}(\Omega) \hookrightarrow H_k \hookrightarrow
  W_k$. In consequence a $W_k$- or a $H_k$-solution of (P).1 is also a
  distributional solution of this equation. Hence $(h,k)$ is a
  distributional solution of (P).
\end{proof}

{\bf Remarks} \noindent \\ 
1. The first point in Proposition \ref{pro1} tells that the notions of 
$H$- or $W$-solution are stronger that the notion of renormalized
solution. This fact is coherent with the second point established in Proposition
\ref{pro1}: a $H$- or $W$-solution is a distributional solution if 
$\nu(k)\in L^1(\Omega)$ whereas a renormalized solution is only a priori
a distributional solution if $\nu(k)\in L^{\infty}(\Omega)$ (see
\cite{lewan} p.185). \\ 
2. if we have $k\in H^1_0(\Omega)$ and if $\nu$ satisfies the growth
condition (\ref{cocor1}) then $\nu(k)\in L^1(\Omega)$.

\subsection*{Comparison with 'energy solution'}

We have seen that when $\nu(k)\in L^1(\Omega)$ then any $W$- (or $H$-)
solution is a distributional solution. Moreover the notion of
$W$-solution coincides with the notion of $H$-solution iff
$W_k=H_k$ (see \cite{zhikov}). \\ 

Some sufficient conditions to have this last equality were established
in \cite{zhikov} and in \cite{gal}, but necessary and sufficient conditions
are not known. \\ 

Let us consider the following condition: 
$$
\text{($R$)}  \quad\left\{
	\begin{array}{l}
	  \sqrt{\nu(k)} \in H^1(\Omega) \\ 
	  T_n(k) \in H^1_0(\Omega), \quad \forall n \in \mathbb N
	  \end{array}
	\right.
$$

It was schown in \cite{gal} that the first condition in ($R$) together
with the property $\nu^{-1} \in L^{\infty}(\mathbb R)$ (which is
assumed in ($H_0$)) implies that $W_k=H_k$. \\ 

In \cite{gal} the authors introduced the notion of 'energy
solution'. They impose ($H_0$) as the basic assumption. Then an 
'energy solution' $(u,k)$ for (P) is in fact a $W$-solution which
satisfies ($R$). This implies that $W_k=H_k$. The energy solution is 
also a $H$-solution, and moreover a distributional solution (because the
first assumption in ($R$) implies that $\nu(k)\in L^1(\Omega)$). \\ 

We see then that the notion of 'energy solution' (in the sense of
\cite{gal}) has the advantage of 
unifying various notions by putting us in the situation where 
$\sqrt{\nu(k)}\in H^1(\Omega)$. The disadvantage is that we have to
impose more complicated conditions on the coefficients $a$ and $\nu$,
in order to obtain a solution. In particular in \cite{gal} Theorem 2.1, 
the authors prove the existence of an 'energy solution' under the
assumptions ($H_0$) and ($H_3$) (see below).
$$
\text{($H_3$)}  \quad\left\{
	\begin{array}{l}
	  \nu \in \mathcal{C}^1(\mathbb R^{+}) \\
	  \exists \ C > 0 \text{ and } \gamma > 1/2 \text{ such that: }\\ 
	  |\nu'(s)| \leq C \quad \forall s \in [0,1] \\ 
	  \frac{|\nu'(s)|}{\sqrt{a(s)\nu(s)}} \leq C.s^{-\gamma} \quad
	  \forall s \geq 1.	  
	\end{array}
	\right.
$$
This condition is not verified in the physical situation ($H_p$), but
only in the approximate situation ($H_p'$).
\\ 

In Theorem \ref{theo1} we obtain a $W$-solution under much 
simpler conditions which are satisfied by ($H_p$). This solution is a
distributional solution under an additionnal simple assumption (see Corollary
\ref{cor1}) 
which is again satisfied in ($H_p$). \\ 

Note also that in the first part of Theorem \ref{theo2} we prove that under the
assumptions ($H_0$) and ($H_2$) (which are satisfied in ($H_p$) if 
$a1\nu_2 = a_2\nu1$), the functions $u$ and $k$ are H\"older 
continuous. In particular $\nu(k) \in L^{\infty}$ which implies that 
$W_k = H_k$, and the notions of $H$-solution, $W$-solution,
distributional solution and renormalized solution coincide in this
case. \\    
 
In order to conclude this Appendix we give a last existence result. Let 
($H_4$) be the following condition: 
$$
\text{($H_4$)}  \quad\left\{
	\begin{array}{l}
	  \nu \in \mathcal{C}^1(\mathbb R^{+}) \\ 
	  \exists \ C > 0 \text{ s.t. } \frac{|\nu'(s)|}{\nu(s)} \leq C
	  \quad \forall s \in \mathbb R.
	\end{array}
	\right.
$$
We have: 

\begin{pro} \label{propfin}

Assume that ($H_0$), ($H_1$) and ($H_4$) hold. Then the
$W$-solution given in Theorem \ref{theo1} is an 'energy solution' (in
the sense of \cite{gal}).

\end{pro}

\begin{proof}
We have assumed that ($H_0$), ($H_1$) hold and consequently all the
results presented in the sections 2, 3 and 4 can be recovered. \\ 

Let $(u,k)$ be the $W$-solution given by Theorem \ref{theo1}.  
By using (\ref{reg1}) we see that the second condition in (R) is
satisfied. Nevertheless we cannot directly conclude that 
$\sqrt{\nu(k)} \in H^1(\Omega)$, but we can obtain a new estimate
for the approximating sequence $(u_n,k_n)$. More precisely, we have: 

\begin{equation} 
\|\sqrt{\nu_n(k_n)}\|_{H^1(\Omega)} \leq C. \label{app1}
\end{equation}

In fact, by using the property that $k_n \in H^1_0(\Omega) \cap
L^{\infty}(\Omega)$ together with $\nu \in \mathcal{C}^1(\mathbb R^{+})$ we obtain 
$\nu(k_n) \in H^1(\Omega)\cap L^{\infty}(\Omega)$, with  
$\nabla \nu(k_n)=\nu'(k_n)\nabla k_n$. Recall now that 
$\nu_n(k_n) = T_n(\nu(k_n))$. Hence we have 
$$\nabla \nu_n(k_n)= 1_{\{\nu_n(k_n) < n\}} \nu'(k_n)\nabla k_n.$$
It follows that: 
\begin{eqnarray*} 
\nabla \sqrt{\nu_n(k_n)}&=& 1_{\{ \nu(k_n) < n\}}
\frac{\nu'(k_n)\nabla k_n}{2 \sqrt{\nu_n(k_n)}}  
= 1_{\{\nu_n(k_n) < n\}}
\frac{\nu'(k_n)}{2 \sqrt{\nu_n(k_n) a_n(k_n)}} \sqrt{a_n(k_n)} \nabla k_n
\\ 
&\overset{\text{by }(\ref{1.8'})}{\leq}& C 1_{\{\nu_n(k_n) < n\}}
\frac{\nu'(k_n)}{\nu_n(k_n)} \sqrt{a_n(k_n)} \nabla
k_n = C \frac{\nu'(k_n)}{\nu(k_n)} \sqrt{a_n(k_n)} \nabla
k_n. 
\end{eqnarray*}
 
Hence, by using (\ref{1.8}) we obtain 
$$\|\nabla \sqrt{\nu_n(k_n)}\|_{L^2(\Omega)} \leq C.$$
Moreover $\sqrt{\nu_n(k_n)} = \sqrt{\nu(0)}$ on $\partial
\Omega$ and thus we obtain (\ref{app1}) by using a Poincar\'e
inequality.
\end{proof}  
{\bf Remark} \noindent \\ 
The hypotheses made in Proposition \ref{propfin} are verified
under assumption ($H_p'$). In the hypotheses, we require only very weak
growth condition at infinity for $\nu$. For instance (contrarily to
the result presented in \cite{gal}) the Proposition \ref{propfin}
works if we have: \\ 
$$\nu(s) = \nu_1 + \nu_2 e^{\beta_1 s}, \quad a(s) = a_1 + a_2
e^{\beta_2 s}, \quad \beta_1 \leq \beta_2.$$  

\section*{Appendix II: H\"older continuity and composition}

Let $\Lambda \subset \mathbb R^d$ and $\alpha \in (0,1)$. We recall that the 
space $\mathcal{C}^{0,\alpha}(\Lambda)$ of H\"older continuous (with exponent $\alpha$) functions on $\Lambda$ is defined by:

$$\mathcal{C}^{0,\alpha}(\Lambda) = \big\{ f \ : \ \Lambda \to \mathbb R 
\quad \text{ s.t. } \ \forall x_0 \in \Lambda: \ 
\sup\limits_{x \in \Lambda} \frac{|f(x)-f(x_0)|}{|x-x_0|^{\alpha}} < \infty 
\big\}$$
More generally, for any integer $k$, the space 
$\mathcal{C}^{k,\alpha}(\Lambda)$ 
is the space of those $f \in \mathcal{C}^k(\Lambda)$ whose kth derivative is in 
$\mathcal{C}^{0,\alpha}(\Lambda)$.\\ 

A first elementary result tells that the product of two H\"older continuous functions is an H\"older continuous function. More precisely we have 
(see relation (4.7) in \cite{gil}): 
\begin{lemm} \label{lem1ap2}
Assume that $f_1, f_2 \in \mathcal{C}^{0,\alpha}(\Lambda)$. Then 
$f_1.f_2 \in \mathcal{C}^{0,\alpha}(\Lambda)$
\end{lemm} 

In Section 5 we used a function defined as a composition of two 
H\"older continuous functions. We needed the following result: 

\begin{lemm}
Let $\overline{\Omega}$ be a compact in $\mathbb R^d$ and $\alpha \in (0,1)$. We consider the following three conditions: 
\begin{eqnarray*}
&(A)& \quad \lambda \in \mathcal{C}^1(\mathbb R) \ \text{ and } \ 
k \in \mathcal{C}^{0,\alpha}(\overline{\Omega}) \\ 
&(B)& \quad \lambda \in \mathcal{C}^{0,\alpha}(\mathbb R) \ \text{ and } \ 
k \in \mathcal{C}^{1}(\overline{\Omega}) \\
&(C)& \quad \lambda \in \mathcal{C}^{1,\alpha}(\mathbb R) \ \text{ and } \ 
k \in \mathcal{C}^{1,\alpha}(\overline{\Omega})
\end{eqnarray*}
We have: \\ 

1. Assume that (A) or (B) is satisfyed. Then $\lambda(k) \in \mathcal{C}^{0,\alpha}(\overline{\Omega})$. \\ 
2. Assume that (C) is satisfyed. Then $\lambda(k) \in \mathcal{C}^{1,\alpha}(\overline{\Omega})$.
\end{lemm}
\begin{proof}
1. In this situation we clearly have $\lambda(k) \in \mathcal{C}^0(\overline{\Omega})$ and  
\begin{equation}
M_1:=\sup\limits_{x \in \overline{\Omega}}|k(x)| < \infty. \label{app2.1}
\end{equation}
Let 
$$I(x,x_0):=\frac{|\lambda(k(x))-\lambda(k(x_0))|}{|x-x_0|^{\alpha}}.$$
We want to prove that 
\begin{equation}
\sup\limits_{x,x_0 \in \overline{\Omega}}I(x,x_0) < \infty. \label{app1}
\end{equation}
$\bullet~$ Assume that (A) holds. Then in addition of (\ref{app2.1}) we have: 
$$M_2:= \sup\limits_{t,t_0 \in [-M_1,M_1]} 
\frac{|\lambda(t)-\lambda(t_0)|}{|t-t_0|} < \infty \quad \text{and} \quad 
M_3:= \sup\limits_{x,x_0 \in \overline{\Omega}} 
\frac{|k(x)-k(x_0)|}{|x-x_0|^{\alpha}} < \infty.$$ 
Consequently: 
$$I(x,x_0) \leq M_2 \frac{|k(x)-k(x_0)|}{|x-x_0|^{\alpha}} \leq M_2.M_3$$
Hence (\ref{app1}) is satisfyed. \\ 
$\bullet~$ Assume now that (B) holds. Then in addition of (\ref{app2.1}) we have: 
$$M_4:= \sup\limits_{x,x_0 \in \overline{\Omega}} 
\frac{|k(x)-k(x_0)|}{|x-x_0|} < \infty \quad \text{and} \quad 
M_5:= \sup\limits_{t,t_0 \in [-M_1,M_1]} 
\frac{|\lambda(t)-\lambda(t_0)|}{|t-t_0|^{\alpha}} < \infty.$$
In this situation we can estimate $I(x,x_0)$ as follows: 
$$I(x,x_0) \leq \frac{|\lambda(k(x))-\lambda(k(x_0))|}{|k(x)-k(x_0)|^{\alpha}}.
\frac{|k(x)-k(x_0)|^{\alpha}}{|x-x_0|^\alpha} \leq M_5 M_4^\alpha.$$ 
Hence (\ref{app1}) is again satisfyed. \\ 

2. Assume that (C) holds and let $\mu:=\lambda(k)$. Clearly 
$\mu \in \mathcal{C}^1(\overline{\Omega})$ and 
$\nabla \mu = \lambda'(k)\nabla k$. \\ 
We remark that $\lambda' \in \mathcal{C}^{0,\alpha}(\mathbb R)$ and  
$k \in \mathcal{C}^{1,\alpha}(\overline{\Omega})$. We can then apply the first point of this lemma to obtain: $\lambda'(k) \in \mathcal{C}^{0,\alpha}(\overline{\Omega})
$. Moreover $\nabla k \in (\mathcal{C}^{0,\alpha}(\overline{\Omega}))^d$. Hence the product $\lambda'(k)\nabla k$ is H\"older continuous.
\end{proof}


\begin{thebibliography}{sto}


\bibitem{num2}
C.\ Bernardi, T.C.\ Rebollo, M.G.\ Marmol, R.\ Lewandowski, 
F.\ Murat: 
'A model for two coupled turbulent fluids Part III: Numerical 
approximation by finites elements',
{\it Num. Math.} {\bf 98}, No 1,(2004), 33--66.

\bibitem{BG}
L.\ Boccardo, T.\ Gallou\"et:
'Non-linear elliptic or parabolic equations involving measure date', 
{\it J. Func. Analysis} {\bf 87}, (1989), 149--169.

\bibitem{brezis}
H. Brezis: 
'Analyse fonctionnelle, th\'eorie et applications', 
Masson, Paris Milan Barcelone (1993).

%\bibitem{bresch}
%D.\ Bresch, J.\ Lemoine, F.\ Guillen-Gonzalez:
%'A note on a degenerate elliptic equation with applications for lakes and 
%seas', 
%{\it Elec. Jour. of Diff. Equa.} {\bf 2004}(42) (2004), 1--13.


%\bibitem{lewan}
%F.\ Brossier, R.\ Lewandowski:
%'Impact of the variations of the mixing length in a first 
%order turbulent closure system',  
%{\it M2AN, Math. Model. Numer. Anal.} {\bf 36}(2) (2002).



%\bibitem{clain}
%S.\ Clain, R.\ Touzani:
%Solution of a two-dimensional stationary induction heating problem without 
%boundedness of the coefficient.
%Math. Mod. Meth. Appl. Sci. {\bf 20} (1997), 759--766.

\bibitem{chen}
Y.-Z. Chen, L.-C. Wu:
'Second Order Elliptic Equations and Elliptic Systems',
Translations of Mathematical Monographs, vol. 174, 
American Math. Society, Providence Rhode Island (1998). 

\bibitem{DL}
R.\ Dautray, J-L.\ Lions:
'Analyse math\'ematique et calcul num\'erique pour les 
sciences et techniques, volume 2, l'op\'erateur de Laplace', 
Masson, Paris New York, (1987).

\bibitem{drey}
P.\ Dreyfuss:
'Higher integrabiblity of the gradient in degenerate elliptic
equations', 
accepted for publication in {\it Potential Analysis}.

\bibitem{gal}
T.\ Gallouet, J.\ Lederer, R.\ Lewandowski, F. Murat, L. Tartar:
'On a turbulent system with unbounded eddy viscosities', 
{\it Nonlinear Analysis} {\bf 52} (2003), 1051--1068.

\bibitem{gil}
D.\ Gilbarg, N.\ S.\ Trudinger:
'Elliptic Partial Differential Equations of Second Order', 
Springer-Verlag, Berlin Heidelberg New York, (1977).

\bibitem{lady}
O.\ A.\ Ladyzenskaja, N.\ N.\ Ural'ceva: 
'Equations aux D\'eriv\'ees Partielles de Type Elliptiques', Dunod,
Paris, 1968.

\bibitem{lewan}
R.\ Lewandowski:
'Analyse math\'ematique et oc\'eanographie', 
Masson Paris, (1997).

\bibitem{lewan2}
J.\ Lederer, R.\ Lewandowski:
'A RANS 3D model with unbounded eddy viscosity',  
to be published in {\it Ann. IHP, An. non Lin.} (2006).

%\bibitem{meyers1}
% N.\ G.\ Meyers: 
%'An $L^p$ estimate for the gradient of solutions of second order 
%elliptic divergence equations',  
%{\it Ann.\ Scu.\ Norm.\ Sup.\ Pisa} {\bf 17} (1963), 189--206.

%\bibitem{MS}
%M.\ K.\ V.\ Murthy, G.\ Stampacchia:
%'Boundary value problems for some degenerate-elliptic operators', 
%{\it Ann. Math. Pura. Appl. IV} {\bf 80} (1968), 1--122.

\bibitem{stam}
G.\ Stampacchia:
'Equations elliptiques du second ordre \`a coefficients 
discontinus', Les presses de l'universit\'e de Montreal, (1966).

\bibitem{cassano}
F. Serra Cassano:
'On the local boundedness of certain solutions for a class of 
degenerate elliptic equations',  
{\it Bollettino U.M.I} (7){\bf 10}-B (1996), 651--680.

\bibitem{zhikov}
V.\ V.\ Zhikov:
'Weighted Sobolev spaces', 
{\it Sbornik Mathematics} {\bf 189}(8) (1998), 1139--1170.

%\bibitem{zhikov2}
%V.\ V.\ Zhikov:
%'To the problem of passage to the limit in divergent nonuniformly 
%elliptic equations', 
%{\it Functional Analysis and Its Applications} {\bf 35}(1) (2001), 19--33.

\end{thebibliography}
\end{document}